\documentclass[journal]{IEEEtran}
\usepackage{ifpdf,cite,epsf,psfrag,epsfig}
\usepackage[cmex10]{amsmath}
\usepackage{mathrsfs, amsthm,amsfonts, latexsym, amssymb,algorithmic,algorithm}
\usepackage[mathscr]{eucal}
\usepackage[caption=false,font=scriptsize]{subfig}
\usepackage{multicol,multirow,color}
\usepackage[breaklinks]{hyperref}
\usepackage[hyphenbreaks]{breakurl}

\newtheorem{lem}{Lemma}
\newtheorem{thm}{Theorem}
\begin{document}
\title{Multi-Area Interchange Scheduling\\under Uncertainty}

\author{Yuting~Ji,~\IEEEmembership{Student~Member,~IEEE,}
        and~Lang~Tong,~\IEEEmembership{Fellow,~IEEE}\thanks {\scriptsize
This work is supported in part by the Department of Energy Office of Electricity Delivery and Energy Reliability Consortium for Electric Reliability Technology Solutions (CERTS) program and the National Science Foundation under Grant CNS-1135844.  Part of this work appeared in \cite{JiTong16PESGM}.

The authors are with the School of Electrical and Computer Engineering, Cornell University, Ithaca, NY 14853, USA (email: yj246@cornell.edu; ltong@ece.cornell.edu).}}
\maketitle
\begin{abstract}
The problem of multi-area interchange scheduling under system uncertainty is considered. A new scheduling technique is proposed for a multi-proxy bus system based on stochastic optimization that captures uncertainty in renewable generation and stochastic load. In particular, the proposed algorithm iteratively optimizes the interface flows using a multidimensional demand and supply functions. Optimality and convergence are guaranteed for both synchronous and asynchronous scheduling under nominal assumptions.
\end{abstract}

\IEEEpeerreviewmaketitle
\section*{Nomenclature}
\addcontentsline{toc}{section}{Nomenclature}
\begin{IEEEdescription}[\IEEEusemathlabelsep\IEEEsetlabelwidth{$C_n(\cdot)$}]
\item[$C_n(\cdot)$]  Cost function of generating units in area $n$ with form $C_n(g_n)=\frac{1}{2}g_n^\intercal H_n g_n+l_n^\intercal g_n$, where $H_n$ is positive definite.
\item[$d_n$]  Vector of net load forecast for area $n$.
\item[$d^t_n$]  Vector of net load forecast for area $n$ at time $t$.
\item[$g_n$]   Vector of dispatch for area $n$.
\item[$q$] Vector of all interface flows with predetermined directions where the $i$th element is denoted by $q(i)$.
\item[$q^{(k)}$] Vector of all interface flows at iteration $k$.
\item[$q^{t}$] Vector of all interface flows at time $t$.
\item[$q_n$] Vector of interface flows associated with area $n$ assuming \textit{outbound} directions.
\item[$q(-i)$] Vector after removing the $i$th entry of $q$.
\item[$F_n$] Vector of transmission capacities for area $n$.
\item[$\mathscr{F}_n$] Distribution of the net load $d_n$, \textit{i.e.}, $d_n\sim \mathscr{F}_n$.
\item[$\mathscr{F}^t_n$] Distribution of the net load $d^t_n$, \textit{i.e.}, $d^t_n\sim \mathscr{F}^t_n$.
\item[$\mathscr{G}_n$] Set of generation constraints for area $n$.
\item[$I$] Number of separate external interfaces.
\item[$N$] Number of independently operated areas.
\item[$Q$] Vector of interface capacities.
\item[$A_n$]  Shift factors of internal buses to lines in area $n$.
\item[$B_n$] Shift factors of proxy buses to lines in area $n$.
\item[$\lambda_n$] Shadow price for area $n$'s power balance constraint.
\item[$\mu_n$] Shadow prices for  area $n$'s transmission constraints.
\item[$\nu_i$] Shadow price for the $i$th interface constraint.
\item[$\pi_n$] Vector of locational marginal prices (LMPs) at proxy buses of area $n$.
\end{IEEEdescription}
\section{Introduction}
The US electric system is partitioned into balancing authority areas connected by tie lines. The interactions among different balancing authorities are carried out through interchange scheduling. Since each balancing authority area, in general, has more than one neighbors, the operator has to specify how much power should flow across each separate scheduling interface. The scheduling of tie line flows is critical to the overall system performance and market efficiency given the fact that the interchange amounts to a significant fraction of net energy. For example, the net flow over external ties accounts for 16.5\% of the net energy for load in ISO New England in 2015\cite{isone_stat}.

In the current external transaction system, operators use two sets of market-based offers to determine external interface schedules: (1) market participants' external transaction offers to buy and sell across an interface, cleared against (2) the real-time generation supply curves in each control area. Inter-regional energy trades are driven by the interface price differentials, facilitated by market participants and administrated by operators. Intuitively, operators can use external tie lines to enhance market efficiencies and operational benefit by transferring excess power from low cost to high cost areas. In practice, however, the state of the art scheduling techniques are not efficient with two common symptoms:  (1) the under utilization of tie lines in transporting power from low cost to high cost areas, and (2) the presence of counter-intuitive flows from high cost to low cost areas. The economic loss from inefficient interchange schedules is substantial, estimated at the level of \$784 million for the New York and New England customers from 2006 to 2010 \cite{IRIS}.

One of the main causes of inefficient interchange scheduling is the latency between scheduling time and physical delivery. Since market participants submit bids and offers for the external transactions at least 75 minutes in advance, the information used for interchange scheduling may not reflect the actual system conditions at the time of physical delivery. With the increasing level of renewable integration, this situation is likely to be exacerbated.

\subsection{Related Work}
Techniques aimed at improving interchange efficiency can be classified into two categories. The first aims to optimize the overall interconnected system in a decentralized fashion. In particular, the optimal interchange schedules are obtained from the multi-area optimal power flow (OPF) problem \cite{KimBaldick97TPS, ConejoAguado98TPS_MOPF, ZhaoLitvinovZheng14TPS,BaldickChatterjee14COR,GuoEtal15PESGM,AhmadiConejoCherkaoui13TPS,LiEtal15TPS}.
Among existing prior work, the authors of \cite{AhmadiConejoCherkaoui13TPS,LiEtal15TPS,DoostizadehEtal16ECM_MAreaParallel} consider the multi-area economic dispatch under wind uncertainty. In \cite{AhmadiConejoCherkaoui13TPS}, a two-stage stochastic market clearing model is formulated for the multi-area energy and reserve dispatch problem whose solution is obtained based on scenario enumerations. In \cite{LiEtal15TPS}, the day-ahead tie-flow scheduling is formulated as a two-stage adaptive robust optimization minimizing the cost of the worst-case wind production. In \cite{DoostizadehEtal16ECM_MAreaParallel}, an adjustable interval robust scheduling of wind
power for day-ahead multi-area energy and reserve market clearing is proposed. The uncertainty of wind farms is represented by predefined intervals and the clearing model is formulated as a mixed integer quadratic
programming problem. For deterministic multi-area economic dispatch approaches, see \cite{KimBaldick97TPS, ConejoAguado98TPS_MOPF, ZhaoLitvinovZheng14TPS,BaldickChatterjee14COR,GuoEtal15PESGM} and reference therein.

The main issue of this category approach is the elimination of arbitrage opportunities for the external market participants. Since operators cannot trade with each other directly, market participants facilitate trades between control areas. The multi-area economic dispatch approach thus cannot be implemented under the current regulation.

The second category includes the current industrial practices based on the so-called proxy bus approximation \cite{IRIS,ChenThorpMount04HICSS,JiZhengTong16TPS}. The proxy bus is a trading location at which market participants can buy and sell electricity. In \cite{ChenThorpMount04HICSS}, a coordinated interchange scheduling scheme is proposed for the co-optimization of energy and ancillary services. The proposal of coordinated transaction scheduling (CTS) in \cite{IRIS} is a state-of-the-art scheduling technique based on the economic argument using supply and demand functions  exchanged by the neighboring operators.  When there is  only a single interface in a two-area system, such functions can be succinctly characterized, and the exchange is only made once; the need of iterations among operators is eliminated.  Built upon the idea of CTS, a stochastic CTS for the two-area single-interface scheduling problem is proposed in \cite{JiZhengTong16TPS}.

A shortcoming of existing techniques based on proxy bus approximations is the difficulty of generalizing it for multi-area interconnected systems where multiple scheduling interfaces have to be optimized simultaneously. The challenge arises from the fact that the interfaces cannot be succinctly characterized by a pair of expected demand and supply functions --- an essential property underlying the approach in \cite{JiZhengTong16TPS}  for the single interface scheduling. When multiple interfaces are involved, the simple idea of equating expected demand and supply functions is not applicable and there is no simple notion that the intersection of demand and supply curves gives the social welfare optimizing interchange.

\subsection{Summary of Contributions}
The main contribution of this paper is twofold.  First, we generalize the single interface scheduling problem to the case involving a network of operating areas, each having multiple interfaces with its neighbors.

We consider two types of interface scheduling: (i) the synchronous scheduling and (ii) the asynchronous scheduling. The former requires all areas operated under the same scheduling clock whereas the latter allows every pair of operating areas setting their interfaces independently of others.  To our best knowledge, there is no existing results on this problem in the open literature.

Second, we present a new scheduling technique based on the classical idea of coordinate descent method. The main idea is to iterate over all interfaces, one at a time, minimizing the overall system cost under uncertainty. In each iteration, a single interface optimization problem is solved by projecting the current solution to a particular coordinate representing a particular interface flow. The optimality and convergence are guaranteed under nominal assumptions for both synchronous and asynchronous algorithms.

\section{Multi-Proxy Bus Representation}
A proxy bus is a physical or virtual location where market participants can trade energy between areas. For interchange scheduling, a proxy bus is the location to which generation in the neighboring area is \textit{assumed} be dispatched up and down in response to the change of interchange schedule. In this paper, we focus on the multi-proxy bus system where more than one proxy buses are used to represent the transmission network of the adjacent systems.
\begin{figure}\centering
\includegraphics[width=.4\textwidth]{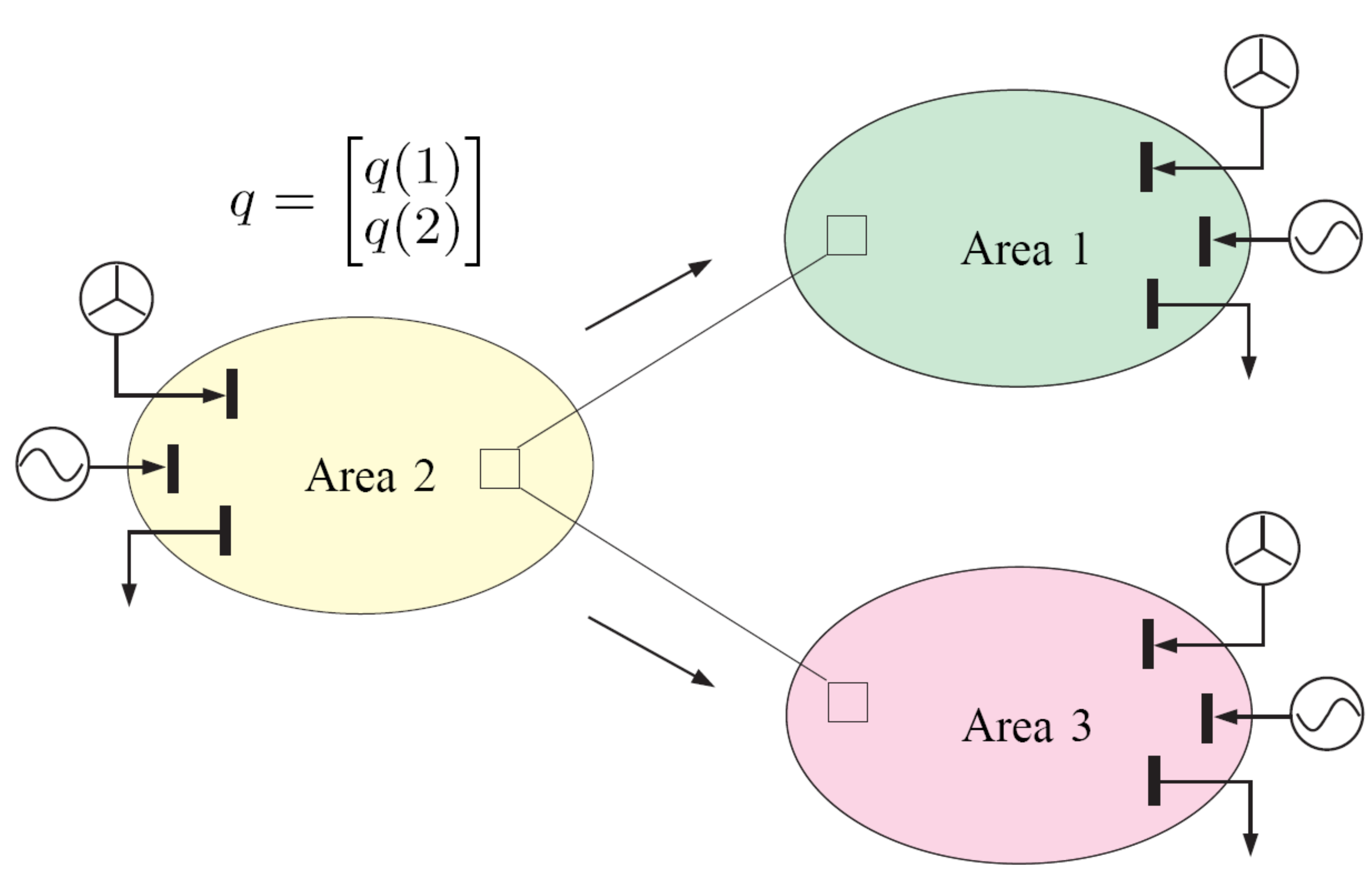}
\caption{A multi-proxy bus representation.}\label{fig:multiproxy}\end{figure}

We use an example to illustrate the multi-proxy bus system. Fig. \ref{fig:multiproxy} depicts a multi-proxy bus representation\footnote{In this system, only the operator of area 2 uses a multi-proxy bus system because it has two separate interfaces. For area 1 and area 3, the approximation is called the single-proxy bus system.} for a 3-area 2-interface interconnected system where each square box indicates a proxy bus. The interchange vector $q$ consists of two interface flows, $q(1)$ and $q(2)$, with fixed directions indicated by the arrows in Fig.~\ref{fig:multiproxy}. Area $n$ maintains a regional interchange vector $q_n$ that describes its own interface flows. For convenience, we assume the direction of each interface flow in the regional interchange vector being \textit{outbound}. In this example, the interchange vector for area 1 is $q_1=-q(1)$, area 2 $q_2=(q(1), q(2))$ and area 3 $q_3=-q(1)$.

The interchange scheduling problem is to optimize the interchange vector $q$ for which the expected overall operation cost is minimized \textit{under the proxy bus model}. This should be distinguished with the problem of multi-area economic dispatch in which the optimal tie line flows and regional generation dispatch are optimized \textit{without network approximation}.

\section{Single Interface Scheduling}
In this section, we consider the single interface scheduling between area 1 and area 2. If the two areas are also physically connected with other areas, we assume that the tie line flows on those areas are fixed, thus not part of the decision process.

The single interface scheduling is a two-stage stochastic optimization: the first stage optimizes the interface flow $q_{12}$ from area 1 to area 2 as in (\ref{opt:s1}); the second stage dispatches internal resources $g^*_1(q_{12},d_1)$ and $g^*_2(q_{12},d_2)$ in the least cost manner to balance the interchange $q_{12}$ and internal net loads $d_1$ and $d_2$ as in (\ref{opt:d3}-\ref{opt:d4}).

The first stage optimization is given by
\begin{equation}\label{opt:s1}
\underset{q_{12}\leq Q_{12}}\min\; \mathbb{E}_{d_1, d_2}\left[C_1\left(g^*_1(q_{12},d_1)
\right)+C_2\left(g_2^*(q_{12},d_2)\right)\right]
\end{equation}
where the expectation is taking over all randomness of internal net loads $d_1\sim\mathscr{F}_1$ and $d_2\sim\mathscr{F}_2$, and $Q_{12}$ is the transfer capacity of the interface between area 1 and area 2.

Given the schedule $q_{12}$, and the realizations $d_1$ and $d_2$, the optimal dispatch $g^*_1(q_{12},d_1)$ and $g^*_2(q_{12},d_2)$ are obtained from the second stage problems:
\begin{equation}\label{opt:d3}
\setlength\arraycolsep{2pt}\begin{array}{r l l}
\underset{g_1\in\mathcal{G}_1}\min &C_1(g_1)&\\
\text{subject to} &\mathbf{1}^\intercal (d_1-g_1)+q_{12}=0 & (\lambda_1)\\
&A_1(d_1-g_1)+b_1 q_{12} \leq F_1 &(\mu_1)\\
\end{array}
\end{equation}
and
\begin{equation}\label{opt:d4}
\setlength\arraycolsep{2pt}\begin{array}{r l l}
 \underset{g_2\in\mathcal{G}_2}\min &C_2(g_2)&\\
\text{subject to} &\mathbf{1}^\intercal (d_2-g_2)- q_{12}=0 & (\lambda_2)\\
&A_2(d_2-g_2)-b_2 q_{12} \leq F_2 &(\mu_2)\\
\end{array}
\end{equation}
where $b_n$, $n\in\{1,2\}$, is the shift factor vector of all transmission lines in area $n$ with respect to the single proxy bus in its neighbor, $\mathbf{1}$ is a vector of ones with the compatible dimension, and the superscript ``$^\intercal$'' denotes the transpose operation.

From the regional dispatch problem, the demand and supply price for the interchange $q_{12}$ can be defined by the cost increase for the individual area with respect to the interchange quantity increase.  By the Envelop Theorem, the LMP $\pi_n(q_{12},d_i)$  at proxy bus in area $n$ is given by
\begin{equation}
\pi_n(q_{12},d_n)=\lambda_n(q_{12},d_n)+b_n^\intercal \mu_n(q_{12},d_n)
\end{equation}
where $\lambda_n(q_{12},d_n)$ and $\mu_n(q_{12},d_n)$ are the Lagrangian multipliers associated with $g^*_n(q_{12},d_n)$ given $q_{12}$ and $d_n$ for area $n$, $n\in\{1,2\}$.

In general, the two-stage stochastic problem is intractable using standard optimization techniques when $d_1$ and $d_2$ follow continuous distributions. Fortunately, for this particular problem (\ref{opt:s1}-\ref{opt:d4}), there is an indirect approach developed in \cite{JiZhengTong16TPS} by solving the following stochastic social welfare maximization problem:

\begin{equation}\label{opt:s2}
\underset{q_{12}\leq Q_{12}}\max \int_0^{q_{12}} \left[\bar{\pi}_2(x)-\bar{\pi}_1(x)\right]dx
\end{equation}
where $\bar{\pi}_n(q_{12})\triangleq\mathbb{E}_{d_n}[\pi_n(q_{12},d_n)]$ is the expected LMP---a function of the net interchange $q_{12}$---at the proxy bus of area $n$, $n\in\{1,2\}$.

The solution of the single interface stochastic optimization problem is given in the following theorem:
\begin{thm}[\hspace{-.05em}\cite{JiZhengTong16TPS}]\label{thm:equivalence}
If  (\ref{opt:d3}) and (\ref{opt:d4}) are not degenerate for all $d_1$ and $d_2$,  then problems (\ref{opt:s1}) and (\ref{opt:s2})  have the same optimizer $q_{12}^*$ satisfying
\[\bar{\pi}_1(q^*_{12})=\bar{\pi}_2(q^*_{12})\] if $q_{12}^*<Q_{12}$ and $Q_{12}$ otherwise.
\end{thm}

Theorem \ref{thm:equivalence} generalizes the tie optimization solution in \cite{IRIS} to the stochastic setting.  Its significance lies in that the optimal interchange is the intersection of expected demand and supply functions rather than the expectation of the intersections of demand and supply curves or the intersection of demand and supply functions using expected generations and demands.

In the following, we generalize the single interface scheduling algorithm to the multi-area system setting where multiple interface flows are involved. In particular, the proposed interface-by-interface scheduling (IBIS) algorithm is specialized for the synchronous scenario in Section \ref{sec:syn} and the asynchronous scenario in Section \ref{sec:asyn}.
\section{Synchronous Interchange Scheduling}\label{sec:syn}
The first scenario we consider is the synchronous interchange scheduling in which all operators in the interconnected system have a unified timetable, \textit{i.e.}, all interface flows are optimized simultaneously at each scheduling time.
\subsection{Problem Formulation}
Consider an interconnected system with $N$ independently operated areas (of an arbitrary network topology) and $I$ separate scheduling interfaces. The multi-area synchronous interchange scheduling problem, analogous to the single interface scheduling problem,  is also a two-stage stochastic optimization: the first stage is to set the values of all interface flows by minimizing the expected overall cost; the second stage is to minimize the cost of individual areas given the fixed interchange and realized random generation and demand.

As a generalization of (\ref{opt:s1}) for the single interface scheduling problem, the first stage optimization for the multi-area system is given by
\begin{equation}\label{p:d1}
\underset{q\leq Q}\min\quad \bar{C}(q)=\sum_{n=1}^N \mathbb{E}_{d_n}\left[C_n\left(g_n^*(q_n,d_n)\right)\right]
\end{equation}
where $q$ is a real vector in dimension $I$, $\bar{C}(q)$ is the expected overall system cost, and $g_n^*(q_n,d_n)$ is the optimal regional dispatch in area $n$, given the interchange level $q_n$ over the interfaces associated with area $n$ and the realized net load $d_n$.

In the second stage, each operator dispatches the internal resource to meet the interchange schedule $q_n$ and the internal net load $d_n$ in the least cost manner subject to the generation and transmission constraints. The optimization problem for area $n$, $n=1,2, \cdots, N$, is specified as
\begin{equation}\label{p:d2}\setlength\arraycolsep{2pt}\begin{array}{r l l}
\underset{g_n\in\mathscr{G}_n}\min &C_n\left(g_n\right)&\\
\text{subject to} &\mathbf{1}^\intercal (d_n-g_n)+\mathbf{1}^\intercal q_n=0 & (\lambda_n)\\
&A_n(d_n-g_n)+B_n q_n \leq F_n. &(\mu_n)
\end{array}\end{equation}

Given the first stage decision $q_n$ and the realization of net load $d_n$, the second stage problems are naturally decoupled and can be solved by their own operators. The LMP vector $\pi_n$ at the proxy buses for area $n$ is calculated from the Lagrangian multipliers of (\ref{p:d2})
\begin{equation}\label{def:pi}
\pi_n(q_n,d_n)=\mathbf{1}\lambda_n^*(q_n,d_n)+B_n^\intercal \mu_n^*(q_n,d_n)
\end{equation}
where $\lambda_n^*(q_n,d_n)$ and  $\mu_n^*(q_n,d_n)$ are functions of the realization  $d_n$ and the first stage decision $q_n$.

The expected multi-dimensional LMP function $\bar{\pi}_n(q_n)$ for area $n$ is therefore defined as
\begin{equation}\label{def:bar(pi)}
\bar{\pi}_n(q_n)=\mathbb{E}_{d_n}[\mathbf{1}\lambda_n^*(q_n,d_n)+B_n^\intercal \mu_n^*(q_n,d_n)]
\end{equation}
where the expectation is taking over all randomness of the net load $d_n\sim\mathscr{F}_n$.
\subsection{Interface-by-Interface Scheduling}
The idea of the proposed scheduling algorithm is to iteratively optimize the interchange vector, one interface at a time, until the termination criterion satisfied. Specifically, at iteration $k$, the $i$th interface flow is given by
\begin{equation}\label{p:s_cmin}
q^{(k)}(i)=\underset{q(i)\leq Q(i)}{\arg\min}\;\bar{C}\left(q(i),q^{(k)}(-i)\right)
\end{equation}
where
\begin{equation}\label{def:qk-i}
q^{(k)}(-i)\triangleq(q^{(k)}(1),\cdots,q^{(k)}(i-1),q^{(k-1)}(i+1),q^{(k-1)}(I)).
\end{equation}

By Theorem \ref{thm:equivalence}, the optimal solution $q^{(k)}(i)$  can be obtained by searching the intersection of the expected supply and demand function defined in (\ref{def:bar(pi)}) for $q(i)$ with fixed $q^{(k)}(-i)$ and check if the interface capacity $Q(i)$ is satisfied.

The detailed synchronous interface-by-interface scheduling (SIBIS) algorithm is given in Algorithm \ref{alg:syn}.

\begin{algorithm}[h]\begin{algorithmic}[1]\caption{Synchronous Interface-by-Interface Scheduling}\label{alg:syn}
\STATE{\textbf{given} a feasible initial point $q^{(0)}$, the expected LMP function $\bar{\pi}_n(q_n)$ for area $n$, $n=1,\cdots, N$, and a tolerance $\epsilon\geq 0$.}
\REPEAT\STATE{$k=k+1$, $q^{(k)}=q^{(k-1)}$.\FOR{$i=1, 2, \cdots, I$}\STATE{Obtain $q^{(k)}(i)$ in (\ref{p:s_cmin}) by intersecting the expected supply and demand functions.}\ENDFOR}
\UNTIL{$\|q^{(k)}-q^{(k-1)}\|_2\leq \epsilon$.}
\end{algorithmic}\end{algorithm}

In practice, a positive value is chosen for $\epsilon$ to ensure a finite termination of SIBIS. When $\epsilon$ is set to zero, the optimality and convergence behavior of SIBIS can be proved (the proof of Theorem \ref{thm:convergence} is given in Appendix \ref{app:thm1}).

\begin{thm}\label{thm:convergence}
Let $\{q^{(k)}\}_{k=0}^\infty$ be the sequence generated by Algorithm \ref{alg:syn} with $\epsilon=0$. Then, every limit point of $\{q^{(k)}\}_{k=0}^\infty$ is optimal to (\ref{p:d1}).
\end{thm}

It should be noted that SIBIS is a form of cyclic coordinate descent method in which one cyclically iterates through the directions, one at a time, minimizing the objective function with respect to each coordinate direction at a time. The early study of the coordinate descent method dates back to 1950s \cite{Hildreth57NRLQ_QP}. The convergence of the method has been extensively studied in the literature \cite{Hildreth57NRLQ_QP,LuoTseng92JOTA_LinearConvergenceCDM,Tseng01_JOTA_BCD} under various assumptions. Given the strict convex assumption of the regional cost function $C_n(g_n)$, the objective function $\bar{C}(q)$ is a continuously differentiable convex function, as shown in the proof of Theorem \ref{thm:convergence}. If $\bar{C}(q)$ has local strict convexity in the feasible region of (\ref{p:d1}), linear rate of convergence can be established as the case in \cite{LuoTseng92JOTA_LinearConvergenceCDM}.

\section{Asynchronous Interchange Scheduling}\label{sec:asyn}
The second scenario we consider is the asynchronous scheduling in which an operator with multiple interfaces determines interface flow one at a time. For such cases, the multi-interface scheduling problem is effectively reduced to a sequential single interface flow optimization.

Mathematically, given the interface flow vector $q^{t-1}(-i)$ at time $t-1$, the schedule of the $i$th interface $q^{t}(i)$ at time $t$  is given by
\begin{equation}\label{p:asyn_single}
q^{(t)}(i)=\underset{q(i)\leq Q(i)}{\arg\min}\sum_{n=1}^N\mathbb{E}_{d^t_n}\left[C_n\left(g^*_n\left(q(i),q^{t-1}(-i),d^t_n\right)\right)\right]
\end{equation}
where the expectation is taking over the randomness of net load $d^t_n$ with respect to the distribution $\mathscr{F}^t_n$ at time $t$, and the regional dispatch $g^*_n\left(q_n),d^t_n\right)$ for area $n$ is the optimal solution to (\ref{p:d2}), given the interchange schedule $q_n$ and the realization of net load $d^t_n$.

Since the interface flow $q(i)$ is the only decision, the objective function  in (\ref{p:asyn_single}) only involves two areas connected by the $i$th interface flow. Therefore, the optimal interface flow  $q^t(i)$  at time $t$ can be obtained by intersecting the expected supply and demand function of $q(i)$ given the distribution $\mathscr{F}^t_n$ of the random net load $d^t_n$. Note that the expected LMP function $\bar{\pi}^t_n(q_n)$, similarly defined in (\ref{def:bar(pi)}), depends on time through the distribution $\mathscr{F}^t_n$.

The distinction between synchronous and asynchronous scheduling lies in the decision at each scheduling time. For synchronous scheduling, the entire interchange vector is optimized via the iterative process given in Algorithm \ref{alg:syn} at each scheduling time $t$. For the asynchronous scheduling, on the other hand, only one element of the interchange vector is optimized at time $t$. Therefore, the solution of the asynchronous scheduling algorithm at time $t$ is suboptimal in terms of minimizing the expected overall system cost.

Below is the description of the asynchronous interface-by-interface scheduling (AIBIS) algorithm where the iterative process is carried out over time which should be distinguished with that in Algorithm \ref{alg:syn}.

\begin{algorithm}[h]\begin{algorithmic}[1]\caption{Asynchronous Interface-by-Interface Scheduling}\label{alg:asyn}
\STATE{\textbf{given} a feasible initial point $q^{(0)}$, the expected LMP function $\bar{\pi}^t_n(q_n)$ for area $n$, $n=1,\cdots,N$, at time $t$, and a termination time $T$.}
\REPEAT\STATE{$t=t+1$, $i=t$ mod $I$.}\STATE{Obtain the optimal solution $q^{t}(i)$ of (\ref{p:asyn_single}) by intersecting the expected supply and demand function.}
\UNTIL{$t=T$.}
\end{algorithmic}\end{algorithm}

We note that if the net load processes $d_n^t$ are independent in time, the optimal interchange depends only on the marginal distribution of the random load at the time of delivery.  If in addition, the process is stationary, \textit{i.e.}, $d_n^t$ is independent and identically distributed (i.i.d.), then the optimal interchange is constant.  In this case, AIBIS is essentially the classical cyclic coordinate decent spread over time.  In comparison with SIBIS, SIBIS achieves optimal interchange for at every time whereas AIBIS achieves the optimality over time.  The following Theorem, whose proof is given in Appedix \ref{app:thm2}, formalize this argument.

\begin{thm}\label{thm:convergence2}
Let $\{q^{t}\}_{t=0}^\infty$ be the sequence generated by Algorithm \ref{alg:asyn} with $T=\infty$. If the net load $d^t_n\overset{i.i.d}\sim\mathscr{F}_n$ for all $n$, then every limit point of  $\{q^{t}\}_{t=0}^\infty$ is optimal to (\ref{p:d1}).
\end{thm}

When $d_n^t$ is not i.i.d., the interchange sequence generated by AIBIS does not converge to that by SIBIS.  The performance (averaged over time) of AIBIS algorithm does not in general converge to that of SIBIS; the lack of synchronization translates to a performance loss.  When the load process is a finite state Markov chain, however, a modification of AIBIS that separately adapts the interchange for different load state will have the same time averaged performance as that of SIBIS.

\section{Evaluation}\label{sec:sim}
In this section, we present numerical results of the proposed scheduling algorithms on the IEEE 118-bus system. The performance of the proposed SIBIS algorithm is compared with the certainty equivalence (CE) technique,\footnote{The CE method also adopts the iterative procedure given in Algorithm \ref{alg:syn}. The only difference lies in the price functions to obtain the single interface flow. Instead of using the expected supply/demand function $\mathbb{E}_{d_n}[\pi_n(q_n,d_n)]$, the CE method uses the supply/demand function $\pi_n(q_n,\mathbb{E}_{d_n}[d_n])$ by substituting the random variable by its the expected value.} which uses the mean value of the random variable to schedule the interchange. The comparison of SIBIS and AIBIS algorithms is then presented for time varying random processes.

\subsection{IEEE 3-Area 118-Bus System}
\begin{figure}
\includegraphics[width=.5\textwidth]{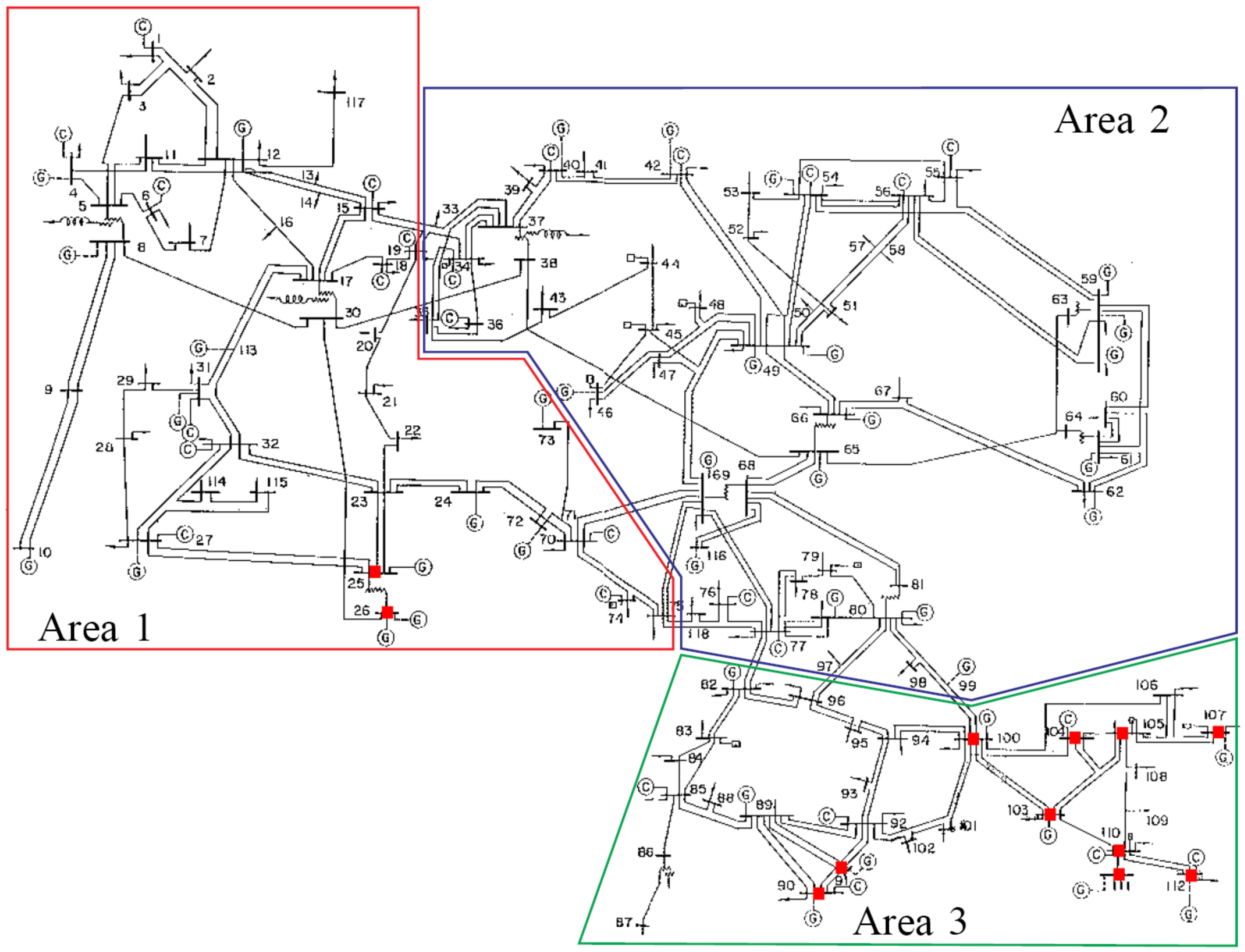}\caption{IEEE 3-area 118-bus system \cite{ieee3area118bus}.}\label{fig:118bus}\vspace{-1em}
\end{figure}
The topology and area partition of the IEEE 118-bus system are given in Fig.~\ref{fig:118bus}. The interchange vector $q$ of this 3-area 2-interface system includes two interface flows where $q(1)$ and $q(2)$ are the flows  to area 2 from area 1 and area 3, respectively. The load profile,  generator capacities, cost functions, and line and bus labels were defined in ``case118'' in \cite{ZimmermanEtal11TPS_Matpower}. We imposed the maximum capacity of $100$ MW on transmission line 8, 126 and 155, and $1200$ MW on the two interfaces. Bus 31, 66 and 92 were selected as proxy buses for area 1, 2, and 3 respectively.

The system uncertainty arose from 12 wind generators (roughly $10\%$ of buses) located at bus 25, 26, 90, 91, 100, 103, 104, 105,  107, 110, 111, and 112, as indicated by red boxes in \ref{fig:118bus}. The selection of these locations were intended to simulate two wind farms; the small one has 2 wind generators in area 1, and the large one has 10 wind generators in area 3 concentrated on a few neighboring buses. All wind generators were assumed to be identical and follow a two-mode Gaussian mixture distribution whose probability density function is given by
\begin{equation}\label{dist:mixgaussian}
f(w)=0.5\frac{\exp\{-\frac{(w-\mu_h)^2}{2\sigma_h^2}\}}{\sqrt{2\pi \sigma_h^2}}+0.5\frac{\exp\{-\frac{(w-\mu_l)^2}{2\sigma_l^2}\}}{\sqrt{2\pi\sigma_l^2}}
\end{equation}
where the Gaussian distribution $\mathcal{N}(\mu_h, \sigma_h^2)$ represents the high wind scenario and $\mathcal{N}(\mu_l, \sigma_l^2)$ the low wind scenario.
\subsection{Expected LMP Computation}
The LMP forecasting technique proposed in \cite{JiTongThomas16TPS} was used to compute the expected LMP. The key idea in \cite{JiTongThomas16TPS} is to solve a multiparametric program to partition the net load space into critical regions\footnote{A critical region is a set of parameters within which the active/inactive status of the inequality constraints of the DC-OPF that determines the LMP are invariant. }.  By the theory of  multiparametric quadratic programming, within each critical region, the LMP is an affine function of the parameter vector. Since the partition of parameter space is independent of realizations of stochastic generation and demand, the computation of critical regions and functions that map the net load to LMP can be obtained ahead of time.  Therefore, the expected LMP $\bar{\pi}_n(q_n)$ can be computed from the conditional distribution of the net load $d_n$ at each level of interchange $q_n$.

Since the optimal interface flow can be obtained by intersecting the expected supply and demand functions, we use the binary search algorithm. Specifically, we did a line search within the interface limits to minimize the expected (sample averaged) price difference.
\subsection{Synchronous Interchange Scheduling}
In this section, we show the optimality and convergence behaviour of the proposed synchronized scheduling algorithm and the two most common symptoms of inefficient schedule.

The distribution (\ref{dist:mixgaussian}) was used as the probabilistic wind production forecast with  parameter values $\mu_h=150$, $\sigma_h=12$, $\mu_l=50$, and $\sigma_l=4$. Note that the difference between the proposed algorithm and the benchmark technique is the use of wind production forecast. The synchronized scheduling algorithm uses the forecasted distribution while the certainty equivalence only uses the mean value $\bar{w}=0.5\mu_h+0.5\mu_l=100$. The initial interchange vector was set at $q^{(0)}=(0,1000)$ for both methods. Termination rule given in Algorithm \ref{alg:syn} was used for both methods with tolerance $\epsilon=0.001$.

\subsubsection{Simulation 1: Tie Line Utilization}  In this simulation, we examine the utilization of tie lines. Ideally, tie lie should be utilized fully to the extent that power flows from a lower price proxy to a high price proxy. All bus loads were set at the default values given in  ``case118'' \cite{ZimmermanEtal11TPS_Matpower}.

From the results shown in Table~\ref{table:iuu}, we observed that the interface flows scheduled by CE and SIBIS had the same directions but different volumes. The CE schedule of transferring $703.32$ MW from area 3 to area 2 resulted in the expected price disparity where the exporting area was \$$31.33$ while the importing area was \$$32.94$. Note that there was adequate transmission capacity, there was still extra transmission capacity and overall cost could have been reduced. This phenomenon is called interface under-utilization, which means transferring more power across the interface can further reduce the overall system cost. The economic benefit from SIBIS can be observed  in the expected overall system cost reduction in the SIBIS schedule. By increasing the interface flow from area 3 to area 2 to $800.16$ MW, the expected supply and demand prices converged to \$$32.44$. Because there was no interface congestion, the expected prices in all three areas converged implying that the SIBIS schedule was efficient.
\begin{table}\renewcommand{\arraystretch}{1.5}\caption{Performance comparison in the interface under-utilization scenario\vspace{-1em}}\label{table:iuu}\begin{center}
\begin{tabular}{c c c c}\hline \hline
&Interface flows (MW)&LMP (\$)&$\mathbb{E}[\text{Cost}]$ (\$)\\ \hline
CE&$\begin{array}{l} q(1)=-68.27\\ q(2)=703.32\end{array}$&$\begin{array}{l} \bar{\pi}_1=32.95\\ \bar{\pi}_2=32.94\\ \bar{\pi}_3=31.33\\ \end{array}$&87114.2\\ \hline
SIBIS&$\begin{array}{l} q(1)=-106.42\\ q(2)=800.16\end{array}$&$\begin{array}{l} \bar{\pi}_1=32.44\\ \bar{\pi}_2=32.44\\ \bar{\pi}_3=32.44\\ \end{array}$&86995.1\\ \hline \hline
\end{tabular}\end{center}\vspace{-1em}\end{table}

\begin{figure}\centering
\subfloat[\vspace{-0.5em}Convergence]{\includegraphics[scale=0.22]{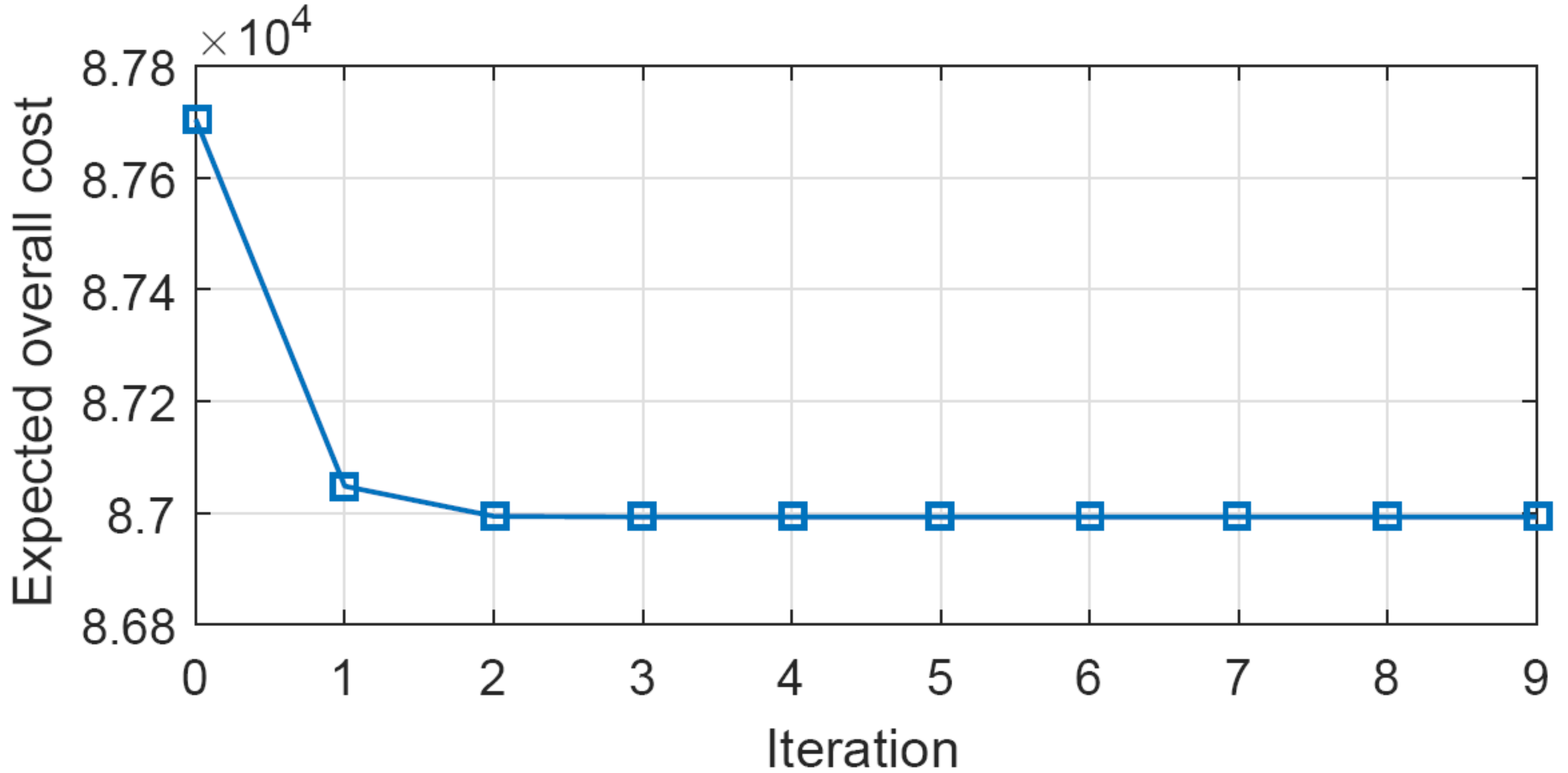}}\\
\subfloat[\vspace{-0.5em}Optimality]{\includegraphics[scale=0.125]{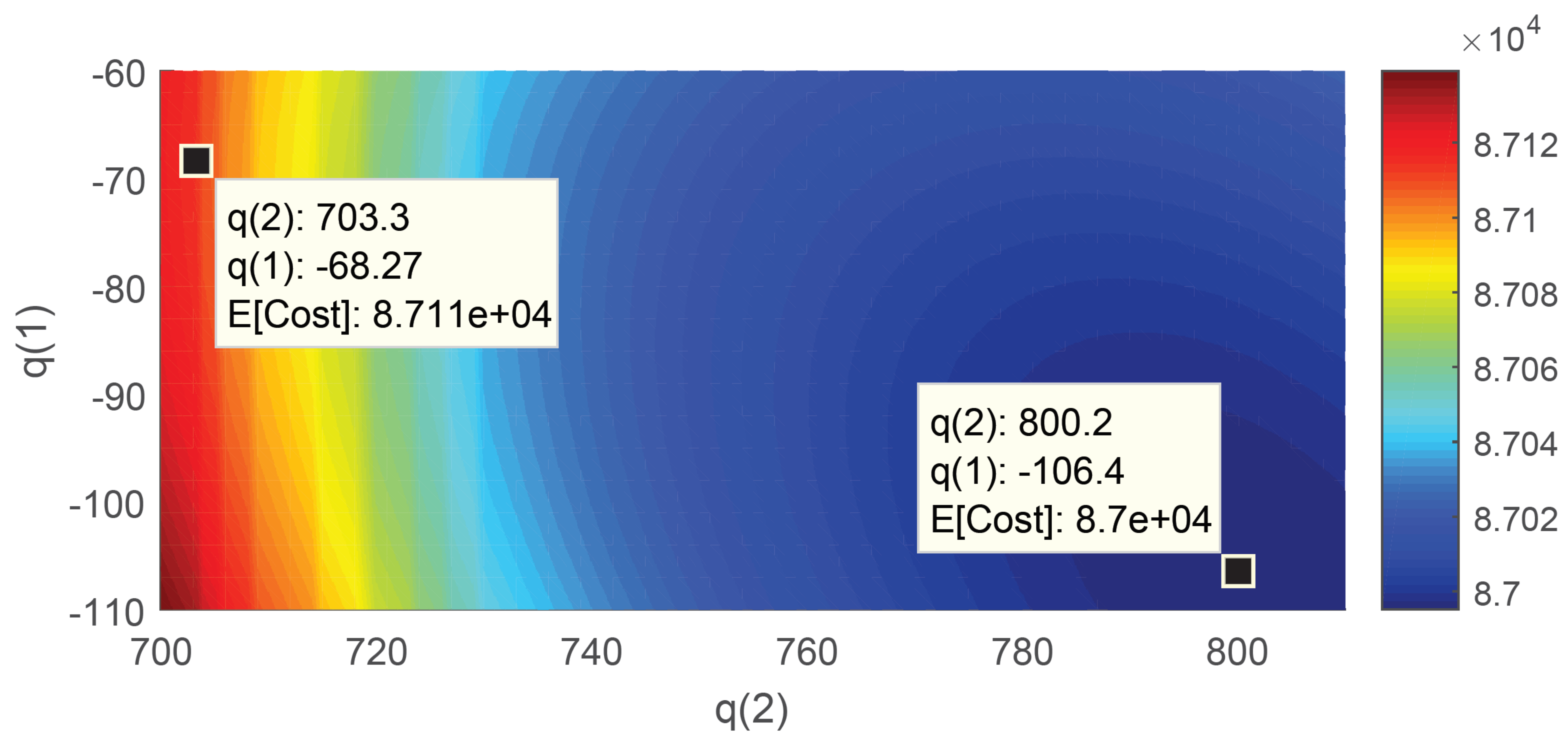}}
\caption{Convergence behavior and optimality of the synchronized scheduling algorithm in interface under-utilization  the scenario.}\label{fig:iuu2}\vspace{-1em}
\end{figure}

The convergence behavior of SIBIS is presented in Fig.~\ref{fig:iuu2}a where the expected overall cost was reasonably close to the optimum after the first two iterations. To demonstrate the optimality of the SIBIS schedule, we computed the expected costs in its neighborhood shown in Fig.~\ref{fig:iuu2}b where the SIBIS schedule is indicated by the cursor at the right bottom and the CE schedule at the left top. Note that the SIBIS schedule is located at the darkest point in this expected cost map which verifies the SIBIS converges to the globally optimal solution.

\subsubsection{Simulation 2: Direction of Tie-line Flow} In this simulation, we examined the direction of power flow on tie lines. In this simulation, the same setting in Simulation 1 was used except that all loads in area 2 were increased by 20\%.

From the results presented in Table~\ref{table:cif}, the interface flows scheduled by the CE method and SIBIS were different in both directions and volumes. Note that the CE method scheduled $56.56$ MW from area 1 to area 2: power flowed from higher- to lower-priced areas. This is economically counter-intuitive. In contrast, the SIBIS schedule resulted in price convergence between proxy buses, thus fundamentally eliminated counter-intuitive flows.
\begin{table}\renewcommand{\arraystretch}{1.5} \caption{Performance comparison in the counter-intuitive flow scenario\vspace{-1em}}\label{table:cif}\begin{center}
\begin{tabular}{c c c c}\hline \hline
&Interface flows&LMP (\$)&$\mathbb{E}[\text{Cost}]$ (\$)\\ \hline
CE&$\begin{array}{l} q(1)=56.56\\ q(2)=790.86\end{array}$&$\begin{array}{l} \bar{\pi}_1=35.07\\ \bar{\pi}_2=34.80\\ \bar{\pi}_3=32.32\\ \end{array}$&101555.2\\ \hline
SIBIS&$\begin{array}{l} q(1)=-2.59\\ q(2)=938.78\end{array}$&$\begin{array}{l} \bar{\pi}_1=34.05\\ \bar{\pi}_2=34.05\\ \bar{\pi}_3=34.05\\ \end{array}$&101364.2\\ \hline \hline
\end{tabular}\end{center}\vspace{-1em}\end{table}

\begin{figure}\centering
\subfloat[Convergence]{\includegraphics[scale=0.22]{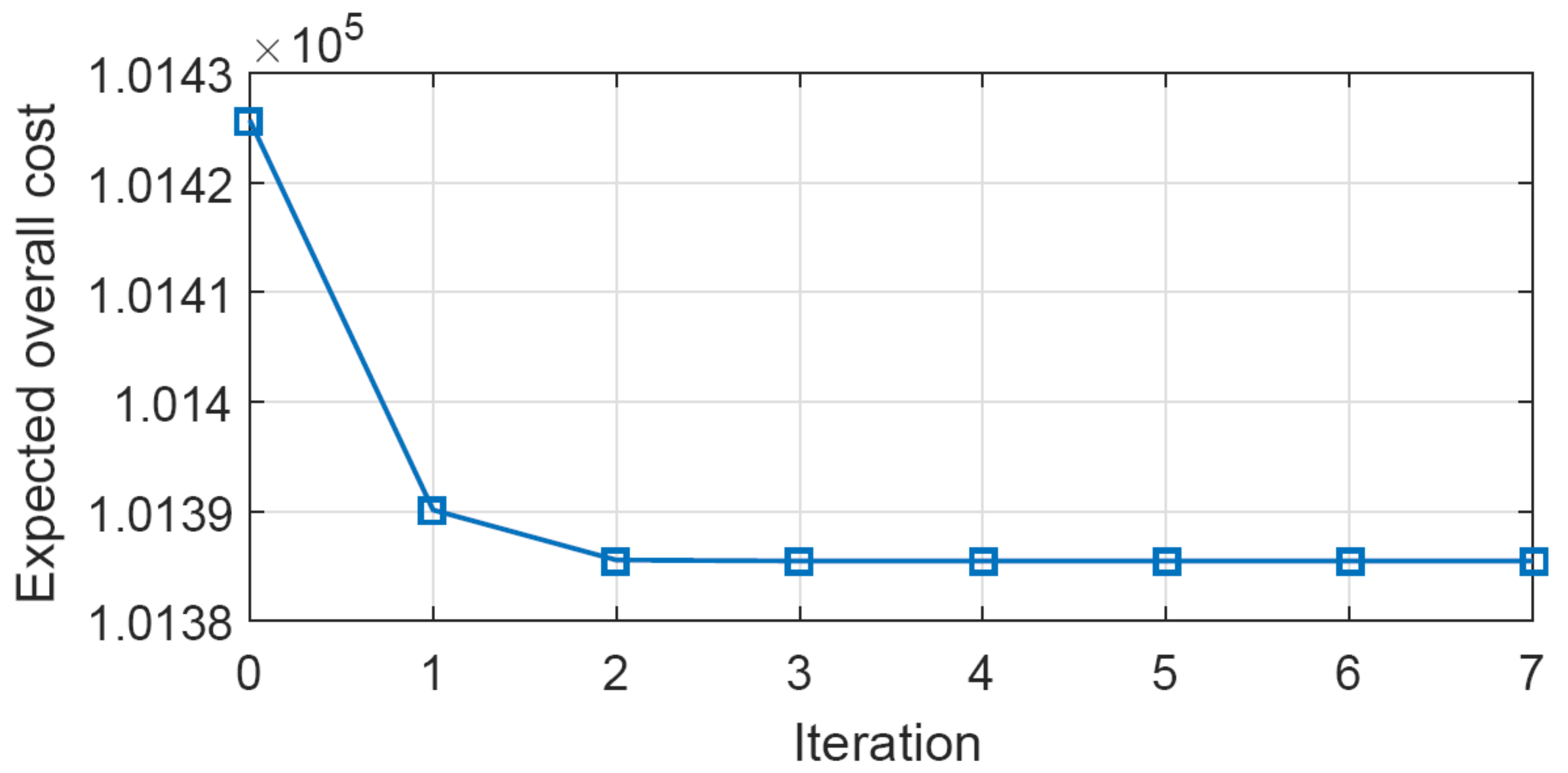}}\\
\subfloat[Optimality]{\includegraphics[scale=0.125]{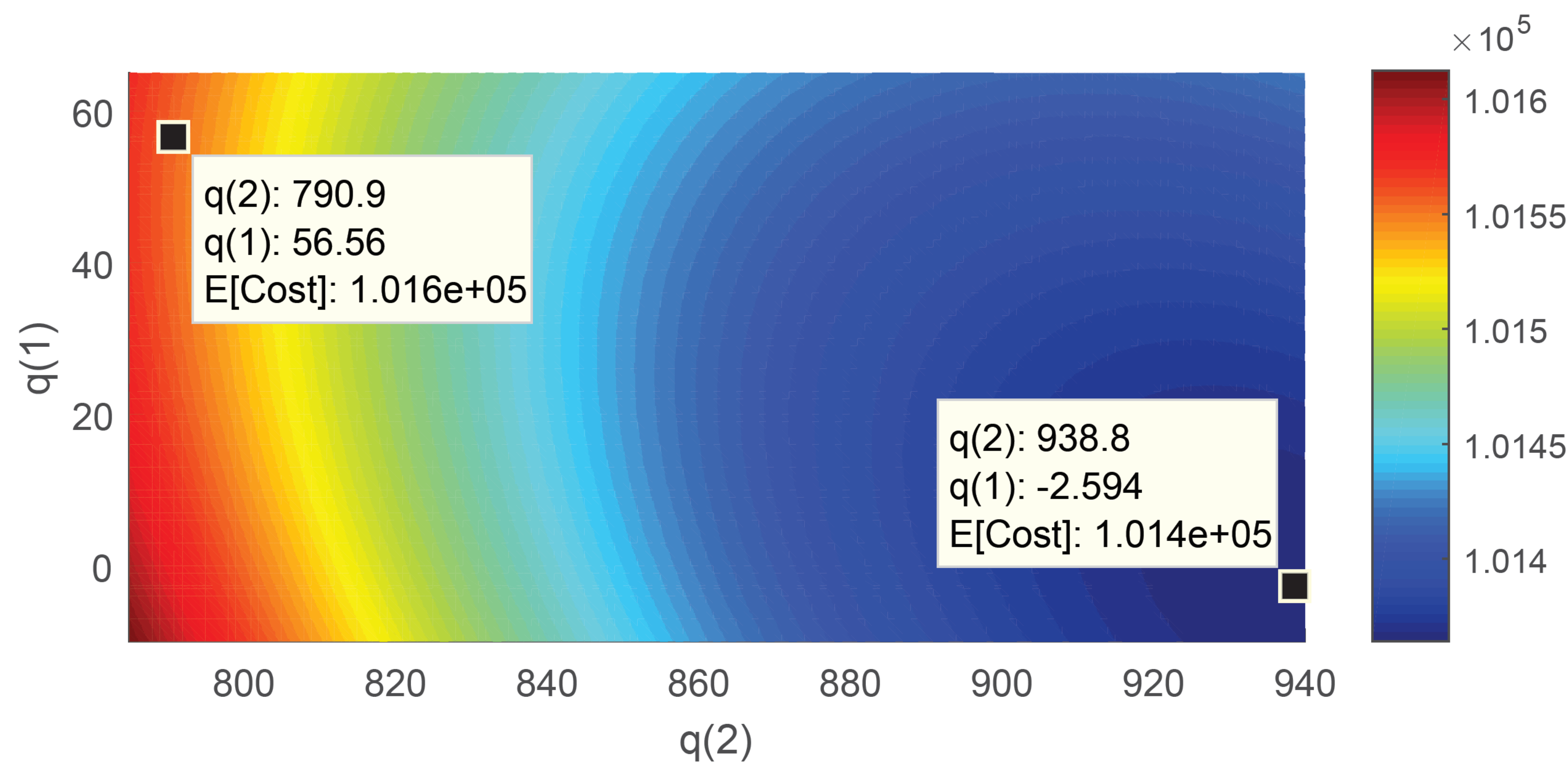}}
\caption{Convergence behavior and optimality of the synchronized scheduling algorithm in the counter-intuitive flow scenario.}\label{fig:cif}\vspace{-1em}
\end{figure}

The convergence behavior was similar to that in Simulation 1 as shown in Fig.~\ref{fig:cif}a. The optimality was also verified.

\subsection{Asynchronous Interchange Scheduling}
We have established theoretical results for AIBIS in Theorem \ref{thm:convergence2} under the independent and stationary assumption for the net load process. In reality, however, the operating condition of the power grid is constantly changing. So we first verify the optimality of AIBIS in the i.i.d. case and compare the performance of SIBIS and AIBIS in a more general setting.

To verify the optimality of AIBIS for the i.i.d. wind process, the generation $w^t(i)$ of wind farm $i$ at time $t$ was assumed to follow $\mathcal{N}(100,10^2)$ for all $i$ and $t$. The scheduling horizon $T$ was set to $20$ for AIBIS and the tolerance $\epsilon$ was set to $0.001$ for SIBIS. The initial interface flows were set at $(500,500)$ for both algorithms.

\begin{figure}\centering
\subfloat[Interface flows]{\includegraphics[scale=0.25]{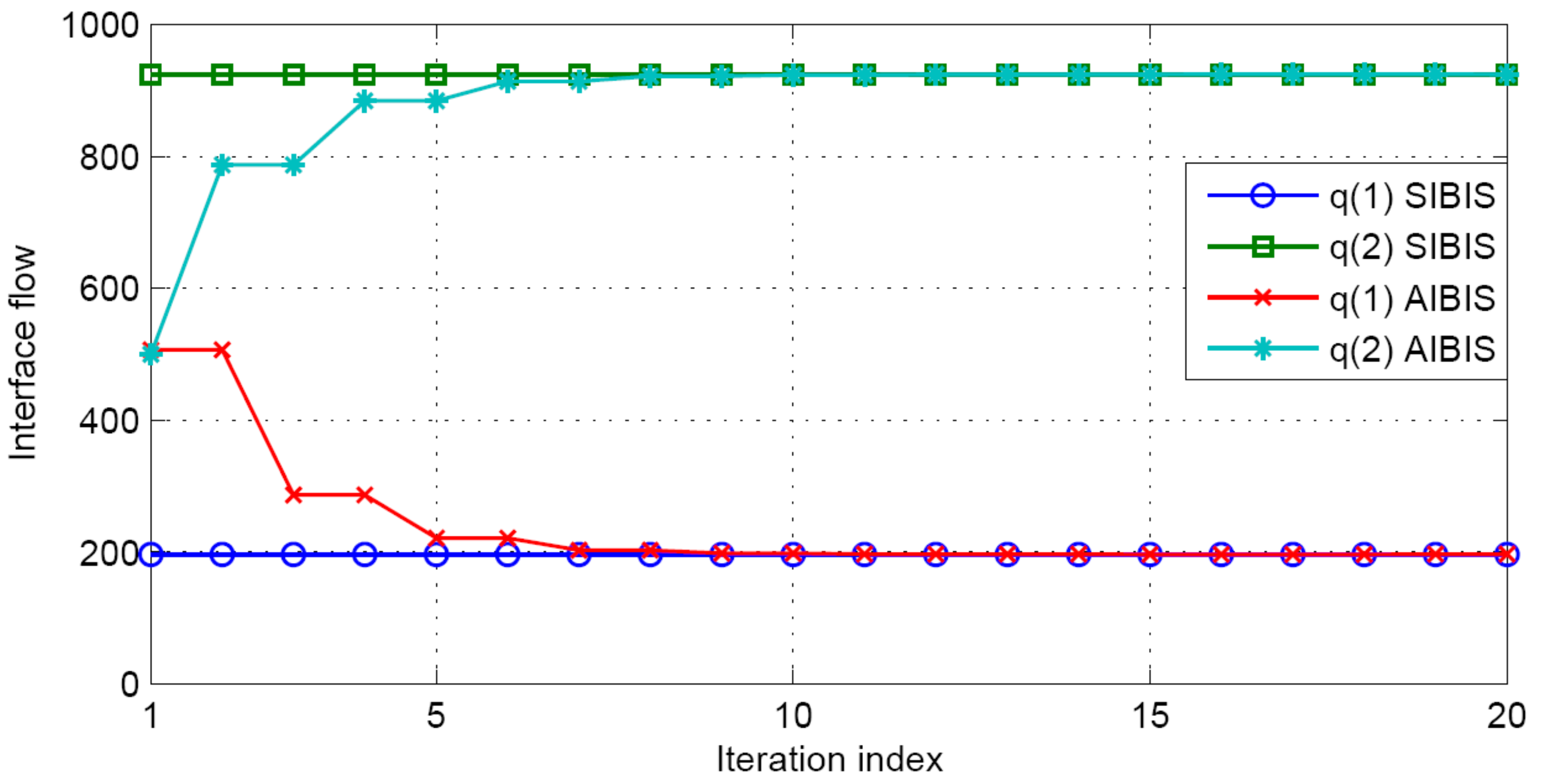}}\\
\vspace{-1em}\subfloat[Expected overall cost]{\includegraphics[scale=0.25]{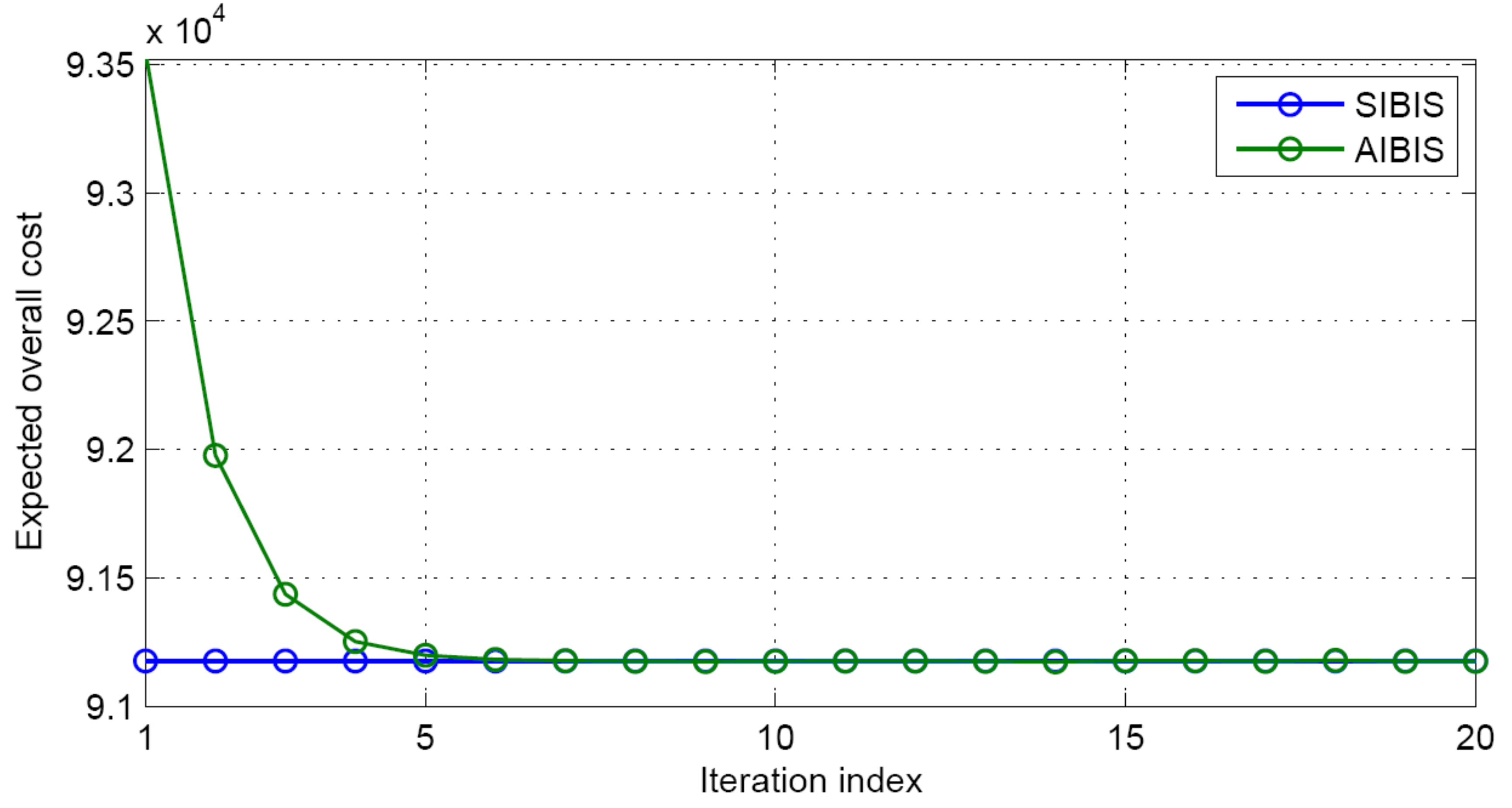}}\\
\caption{Performance comparison for SIBIS and AIBIS with i.i.d. wind generations.}\label{fig:asyn_iid}\vspace{-1em}
\end{figure}

As shown in Fig.~\ref{fig:asyn_iid}, SIBIS achieves optimal interchange for at every time whereas AIBIS achieves the optimality over time. Specifically, the interface flows and the expected cost of SIBIS are optimal and remain constant over time. AIBIS generated a sequence that converged to the optimal interchange and expected cost at time 7. It should be noted that the convergence rate is highly dependent on the initial values. The convergence time of AIBIS is simply the adaptation of the cyclic coordinate descent method.

The behavior of the proposed AIBIS algorithm was then investigated using a time varying process of wind generation. Specifically, we varied the mean value of wind generation starting from $100$ MW and ending at $140$ MW with a constant increment $2$ MW. Except for the varying mean, the rest of the setting remained the same.
\begin{figure}\centering
\subfloat[Interface flows]{\includegraphics[scale=0.25]{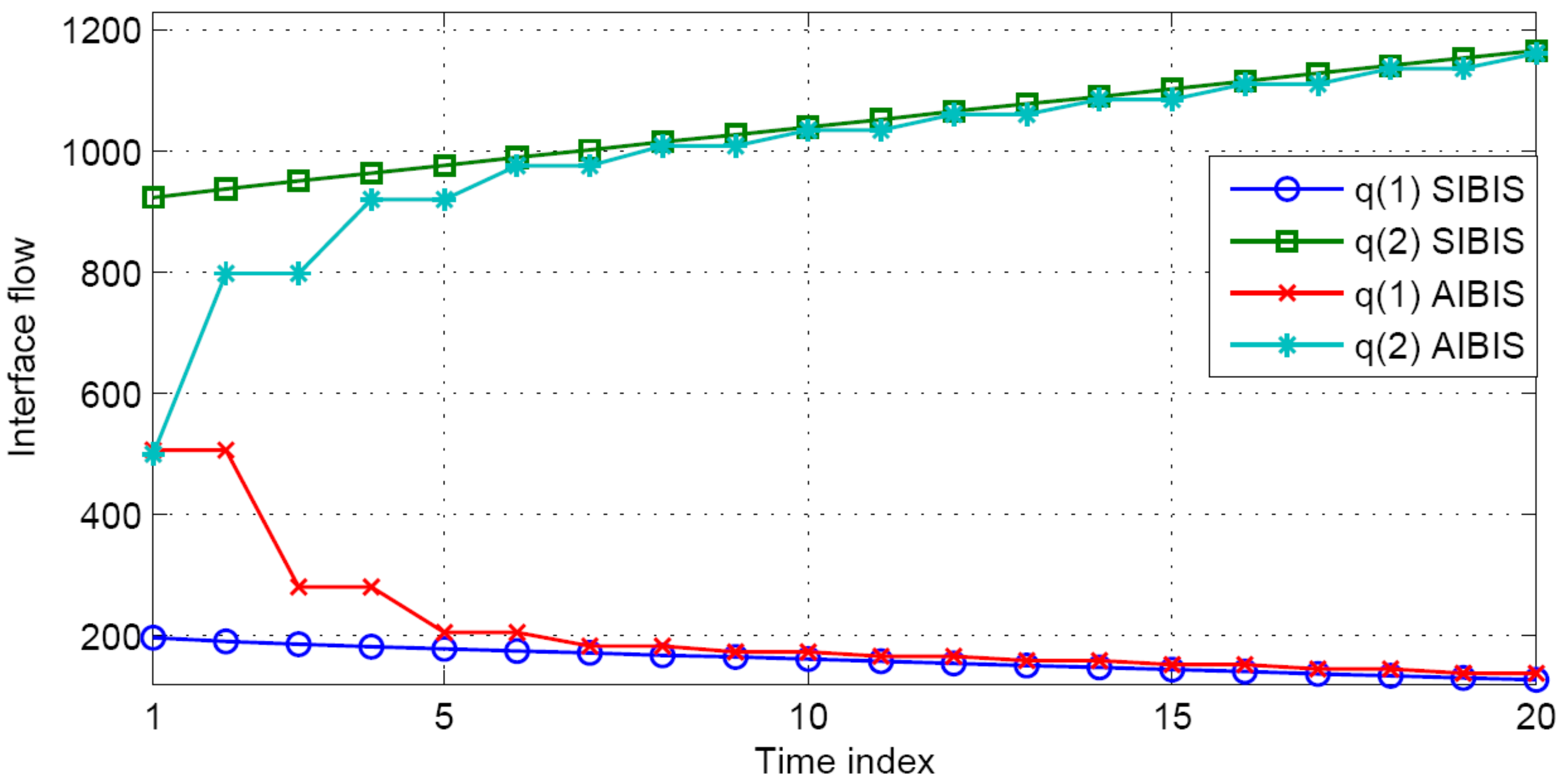}}\\
\vspace{-1em}\subfloat[Expected overall cost]{\includegraphics[scale=0.25]{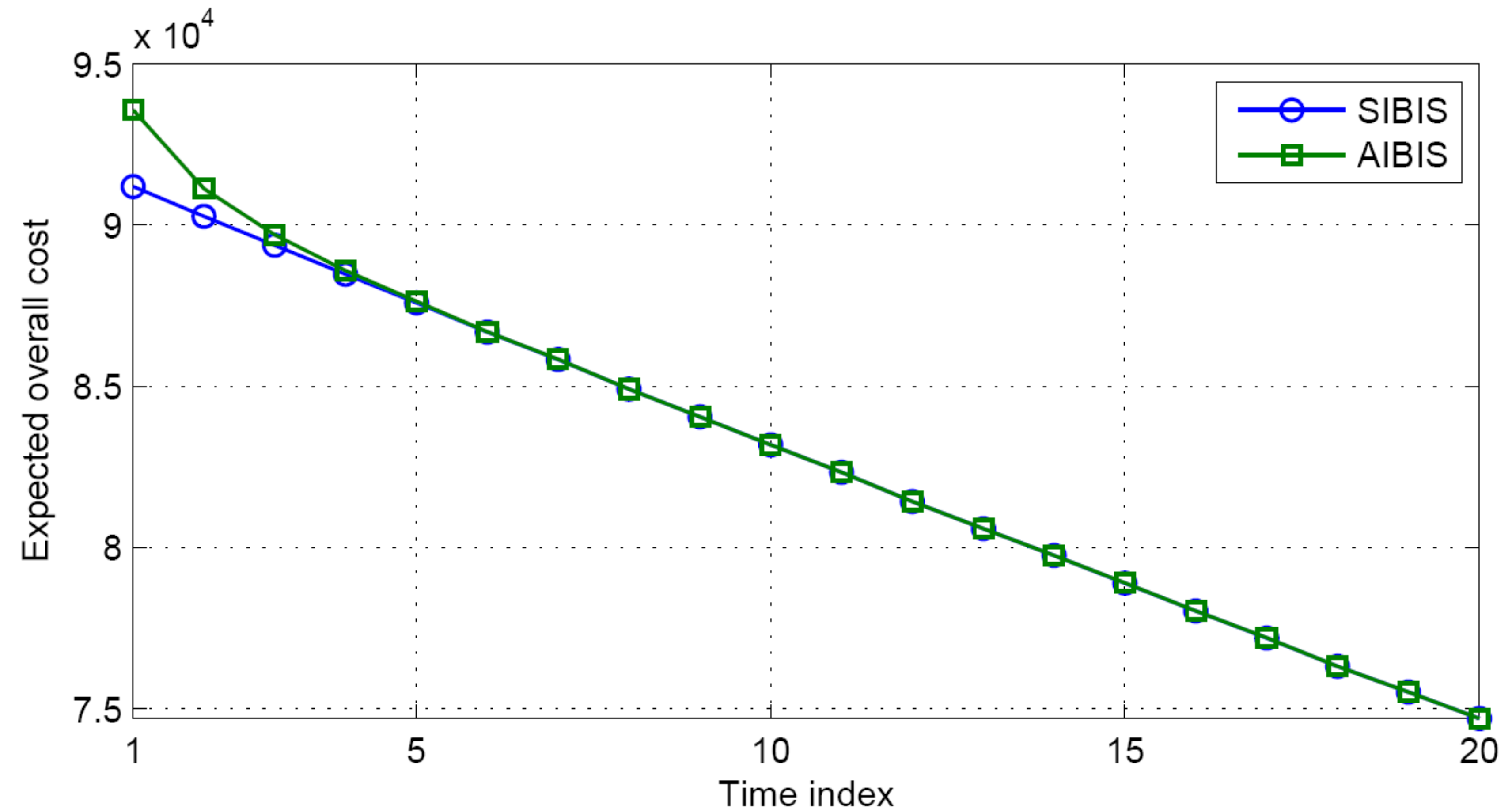}}\\
\caption{Performance comparison for SIBIS and AIBIS with time-varying wind profiles.}\label{fig:asyn}\vspace{-1em}
\end{figure}

From Fig.~\ref{fig:asyn}a, the difference between SIBIS and AIBIS was clearly observed in the scheduled interface flows. Since there were only two interfaces, the interface flows scheduled by AIBIS were alternatively constant during a scheduling time slot. With time increasing, the scheduling difference decreased and so did the expected overall cost shown in Fig.~\ref{fig:asyn}b.

\section{Conclusion}\label{sec:con}
This paper presents a new interchange scheduling technique with the consideration of load and generation uncertainties in the context of multiple proxy bus system. Built upon the idea of coordinate descent method, the interchange vector is iteratively determined, one at a time, to minimize the expected overall system cost.
\appendices
\section{Proof of Theorem \ref{thm:convergence}}\label{app:thm1}
\begin{proof}
We first give the properties, implied by Lemma \ref{lem:mp_affine}, of the objective function $\bar{C}(q)$ in (\ref{p:d1}) and establish the convergence of the algorithm.
\begin{lem}[\hspace{-.05em}\cite{JiZhengTong16TPS}]\label{lem:mp_affine}
If problem (\ref{p:d2}) is neither primal nor dual degenerate for all $d_n$, then the optimal value function $C_n^*(q_n,d_n)\triangleq C_n(g_n^*(q_n,d_n))$ is convex, continuously differentiable and piece-wise quadratic in both $q_n$ and $d_n$.
\end{lem}
Let $q^{(k_i)}=\left(q^{(k)}(i),q^{(k)}(-i)\right)$ where $q^{(k)}(-i)$ is defined in (\ref{def:qk-i}). Using the updating rule (\ref{p:s_cmin}), we have
\begin{equation}\label{eqn:mono}
\bar{C}(q^{(k)}) \leq \bar{C}(q^{(k_{I-1})})\leq\cdots \leq \bar{C}(q^{(k_1)})\leq \bar{C}(q^{(k-1)}), \forall k.
\end{equation}
Let $q^*$ be a limit point of the sequence $\{q^{(k)}\}_{k=1}^\infty$. Note that $q^*\in\mathscr{Q}=\{q|q\leq Q\}$ because $\mathscr{Q}$ is closed. The monotonicity (\ref{eqn:mono}) implies that the sequence $\{\bar{C}(q^{(k)})\}_{k=1}^\infty$ converges to $\bar{C}(q^*)$.

Let $\{q^{(k)}\}_{k\in\mathscr{K}}$, where $\mathscr{K}$ is an index set, be a subsequence of $\{q^{(k)}\}_{k=1}^\infty$ that converges to $q^*$. From the updating rule (\ref{p:s_cmin}) and the monotonic property (\ref{eqn:mono}), for any interface $i$, we have
\[\bar{C}\left(q^{(k)}\right)\leq \bar{C}\left(q^{(k_i)}\right)\leq \bar{C}\left(q(i),q^{(k)}(-i)\right), \forall q(i)\leq Q(i). \]
Since $\bar{C}(q)$ is continuous, implied by Lemma \ref{lem:mp_affine}, taking the limit as $k\in\mathscr{K}$ tends to infinity on both sides, we have
\[\bar{C}(q^*)\leq \bar{C}\left(q(i),q^*(-i)\right), \forall q(i)\leq Q(i)
\]
which means $q^*(i)$ is an optimal solution of the following optimization
 \begin{equation}\label{opt:q(i)}
 \underset{q(i)\leq Q(i)}{\min}\bar{C}\left(q(i),q^*(-i)\right).
 \end{equation}
Therefore, $q^*(i)$ satisfies the Karush-Kuhn-Tucker (KKT) conditions for (\ref{opt:q(i)}), \textit{i.e.}
\setlength\arraycolsep{2pt}\begin{eqnarray}\label{eqn:kkt1}
\nabla_i \bar{C}(q^*(i),q^*(-i))+\lambda(i)(q^*(i)-Q(i))&=&0\\
q^*(i)&\leq& Q(i)\\\label{eqn:kkt2}
\lambda(i)&\geq& 0\label{eqn:kkt3}
\end{eqnarray}
where $\nabla_i \bar{C}(q)$ is the partial derivative with respective to $q(i)$ and $\lambda(i)$ the associated Lagrangian multiplier.

Note that conditions (\ref{eqn:kkt1}-\ref{eqn:kkt3}) hold for all $i$ at $q^*$, \textit{i.e.},
\begin{eqnarray}\label{eqn:kkt21}
\nabla \bar{C}(q^*)+\lambda(q^*-Q)&=&0\\
q^*&\leq& Q\\
\lambda&\geq& 0.\label{eqn:kkt23}
\end{eqnarray}

Since conditions (\ref{eqn:kkt21}-\ref{eqn:kkt23}) are the KKT conditions for (\ref{p:d1}), and $\bar{C}(q)$ is convex by Lemma \ref{lem:mp_affine},  $q^*$ is optimal to (\ref{p:d1}).
\end{proof}

\section{Proof of Theorem \ref{thm:convergence2}}\label{app:thm2}
\begin{proof}
Since $d^t_n\overset{i.i.d}\sim\mathcal{F}_n$ for all area $n$, $q^t(i)$ defined in (\ref{p:asyn_single}) can be obtained by
\begin{equation}\label{eqn:qt(i)}
q^t(i)=\underset{q(i)\leq Q(i)}{\arg\min}\bar{C}\left(q(i),q^{t-1}(-i)\right)
\end{equation}
which implies \[\bar{C}(q^t)\leq \bar{C}(q^{t-1}), \forall t.\]
Let $\tilde{q}$ be a limit point of the sequence $\{q^{t}\}_{t=1}^\infty$. The monotonicity of $\bar{C}(q^{t})$ implies that the sequence $\{\bar{C}(q^{t})\}_{t=1}^\infty$ converges to $\bar{C}(\tilde{q})$.

Let $\{q^{t}\}_{t\in\mathscr{T}}$, where $\mathscr{T}$ is an index set, be a subsequence of $\{q^{t}\}_{t=1}^\infty$ that converges to $\tilde{q}$. From the updating rule (\ref{eqn:qt(i)}) and the monotonicity of $\bar{C}(q^{t})$, we have
\[\bar{C}\left(q^{t}\right)\leq\bar{C}\left(q(i),q^{t-1}(-i)\right), \forall t, \forall i, \forall q(i)\leq Q(i). \]
Since $\bar{C}(q)$ is continuous, implied by Lemma \ref{lem:mp_affine}, taking the limit as $t\in\mathscr{T}$ tends to infinity on both sides, we have
\[\bar{C}(\tilde{q})\leq \bar{C}\left(q(i),\tilde{q}(-i)\right), \forall q(i)\leq Q(i)
\]
which means $\tilde{q}(i)$ is an optimal solution of the following optimization
 \begin{equation}\label{opt:qt(i)}
 \underset{q(i)\leq Q(i)}{\min}\bar{C}\left(q(i),\tilde{q}(-i)\right).
 \end{equation}
Therefore, $\tilde{q}(i)$ satisfies the KKT conditions (\ref{eqn:kkt1}-\ref{eqn:kkt3}) for \ref{opt:qt(i)} at $\tilde{q}$. Since (\ref{eqn:kkt1}-\ref{eqn:kkt3}) hold for all $i$ at $\tilde{q}$, $\tilde{q}$ satisfies the KKT conditions (\ref{eqn:kkt21}-\ref{eqn:kkt23}). By the convexity of $\bar{C}(q)$, the KKT conditions are sufficient and necessary for optimality. Therefore, $\tilde{q}$ is optimal to (\ref{p:d1}).
\end{proof}
\bibliographystyle{IEEEtran}
{\bibliography{reference}}
\end{document}